\documentstyle{amsppt}
\magnification =\magstep1\hsize 6.5truein
\vsize 9truein
\expandafter\redefine\csname logo\string@\endcsname{}
\topmatter
\title The Golden Age of Immersion Theory in Topology: 1959-1973
\endtitle
\author David Spring \endauthor
\TagsOnRight
\endtopmatter
\define\pt{\partial}
\define\Xr{X^{(r)}}
\define\ep{\epsilon}
\document

\flushpar {\bf 1.\,Brief Overview.} In this article \footnote"*"{A preliminary version of this article was presented at the  meeting of the Canadian Mathematical Society, History of Mathematics Session, December 2001.} I briefly review selected contributions to immersion-theoretic topology during the early ``golden'' period, from about 1959-1973, during which time the subject received its initial important developments from leading topologists in many countries. Modern immersion-theoretic topology began with the work of Stephen Smale (1958), (1959) on the classification of immersions of the sphere $S^n$ into Euclidean space $\bold R^q$, $n\ge 2$, $q\ge n+1$. During approximately the next 15 years the methods introduced by Smale were generalized in various ways to analyse and solve an astonishing variety of geometrical and topological problems. In particular, at Leningrad University during the late 1960s and early 1970s, M. Gromov and Y. Eliashberg developed geometrical methods for solving general {\it partial differential relations} in jet spaces where, informally, a partial differential relation of order $r$, $r\ge 1$, is a set of equations (closed relations) or inequalities (open relations) on the partial derivatives of order $\le r$ of a function. For example the immersion relation above studied by Smale is an open 1st order partial differential relation: the function $f\: S^n\to\bold R^q$ is an immersion if in local coordinates $x=(x_1,\dots ,x_n)$ at each point of $S^n$ the first order partial derivatives $\pt_{x_i}f(x)\in\bold R^q$ are linearly independent, hence the Jacobian matrix for $Df(x)$ has maximal rank $(=n)$, which is an open condition on the 1st partial derivatives. Immersion-theoretic topology is a name given to the geometrical methods and to the body of topological results that have been developed over the past 40 years to study the existence and classification, up to homotopy, of global solutions to partial differential relations. There does not appear to be in the literature a general history of immersion-theoretic topology. Though not comprehensive, I highlight topics and issues that reflect my particular interest over a period of many years in solving the $h$-principle (cf. \S1.3) for partial differential relations that occur in various topological contexts. 
    
As indicated above, modern immersion-theoretic topology began with the pioneering work of Smale (1958) on the classification of immersions of the 2-sphere in $\bold R^q$, $q\ge 3$, followed by his generalization Smale (1959) to the classification of immersions of the $n$-sphere in $\bold R^q$, $q\ge n+1\ge 3$. Smale is a co-recipient of the 1966 Fields Medal, for which his pioneering work in immersion theory was included in the citation. For a detailed, comprehensive review of Smale's contributions to differential topology cf. Hirsch (1993). Important earlier works include Whitney (1944), who proved that all smooth $n$-manifolds can be immersed into $\bold R^{2n-1}$, $n\ge 2$, and John Nash (1954) who solved the $C^1$-isometric immersion problem. Nash proved that an abstract continuous Riemannian metric on a smooth manifold $M^n$ can be realized concretely as the induced metric on the manifold obtained through a suitable immersion of class $C^1$ of the manifold into $\bold R^q$, $q\ge n+2$. Some refinements were introduced by Kuiper (1955) who improved these results to $q\ge n+1$. The geometrical methods introduced in Nash (1954) were seen as a curiosity in topological circles of the time until the work of Gromov (1973b) on convex integration theory for which, according to Gromov, Nash (1954) was a precursor (cf. Spring (1993) for further historical remarks). John Nash is a co-recipient of the 1994 Nobel Prize in economics, based on his innovative work in non-cooperative game theory. Recently the life of John Nash was the subject of the Hollywood film ``A Beautiful Mind'', released in 2001. Misha Gromov has received  many awards, including the Wolf Prize (1993) and the Kyoto Prize (2002) for his contributions to mathematics. 

During the late 1960s the story of immersion-theoretic topology shifts to the topology seminar of V. Rokhlin at Leningrad University (now St. Petersburg). Rokhlin is known for his work in ergodic theory and also for his work in the topology of 4-manifolds, in particular for the development of what is known today as the Rokhlin invariant. Rokhlin was interested in the new results on immersion theory that were developed by Smale and other researchers in the West. One of his promising pupils during the mid 1960s was Misha Gromov, whose thesis would generalize Smale's theory of immersions to corresponding weak homotopy equivalence results in the context of very general open relations defined in higher order jet spaces. Towards the end of the 1960s, another promising student, Yasha Eliashberg, began attending Rokhlin's seminar while he was still an undergraduate student (apparently this was not uncommon in Russian circles). Eliashberg and Gromov became friends and they exchanged important ideas about immersion-theoretic topology; this collaboration led to their first joint papers together, on the topic of the removal of singularities, discussed briefly below. Eliashberg's thesis (1972) on the  surgery of singularities was another major advance in this topic as it applies to the simplification of singularities of maps. Shortly afterwards, Gromov published his first paper on convex integration theory (1973b) which took the subject to another conceptual level with applications to diverse new problems in topology and also non-linear partial differential equations. These mathematical developments, from Smale (1959) to Gromov (1973b), serve as natural milestones that enclose the Golden Age of Immersion Theory, the title of this paper. I should add that the subject of immersion-theoretic topology did not come to a close in 1973. Indeed, there have been new results and new perspectives on immersion-theoretic topology in recent years that I will mention briefly throughout this paper and also in a final postscript. 

Historically, I mention also that around 1974 Gromov managed to leave Russia for a position at SUNY Stony Brook. In the early 1980s he moved to Paris, France, and in 1982 he became a permanent member of the Institut des Hautes Etudes Scientifiques (I.H.E.S.). After a series of professional hardships Eliashberg also managed to leave Russia in the late 1980s, later settling at Stanford University. Both Gromov and Eliashberg have remained on the forefront of modern topology, in particular symplectic and contact topology, for the past 10-15 years. Their close collaboration dates back to Rokhlin's seminar in Leningrad. 
\vskip .25cm
\flushpar{\bf 2.\,Smale's Theory.} A smooth map between manifolds, $f\: N^n\to W^{q}$, $q\ge n$, is an {\it immersion} if the map $f$ is of maximal rank ($=n$) at each point of $N^n$. Locally, in suitable curvilinear coordinates near each point of $N^n$ this means that $f$ is modeled on the inclusion map $i\: \bold R^n\to\bold R^q$, whereas globally there is no restriction on the self intersections set of $f$. For example a figure 8 in the plane is an immersed circle (with one double point). An immersion $f\: N^n\to W^q$ is an {\it embedding} if in addition $f\: N^n\to f(N^n)\subset W^q$ is a homeomorphism onto its image. Two immersions $f,g\: N^n\to W^q$ are {\it regularly} homotopic if there is a smooth homotopy of immersions connecting $f$, $g$\,: a smooth map $F\: N^n\times [0,1]\to W^q$, $F_0=f$, $F_1=g$, such that for each $t\in[0,1]$, the map $F_t\: N^n\to W^q$ is an immersion. Thus the tangent map $dF_t\: T(N^n)\to T(W^q)$ is smooth in the variable $t\in[0,1]$\,; in particular, the tangent planes to the immersion $F_t\: N^n\to W^q$ vary continuously in the variable $t\in[0,1]$, an important point when studying examples. The classification problem consists of classifying algebraically, usually in terms of homotopy invariants, immersions of $N^n$ into $W^q$, up to regular homotopy. The first case of interest is the one-dimensional case: immersions of the circle $S^1$ into the plane $\bold R^2$ were classified much earlier by Whitney (1937) in terms of the winding number $\in \bold Z$ of the unit tangent vector to the immersed circle. Of central importance to the proof of this theorem is the fact that one can integrate an appropriate homotopy of tangent vectors to obtain a regular homotopy of immersions of $S^1$. This in turn is related to the classical fact that smooth vector fields can be integrated to obtain integral curves. Attempts to generalize this technique to integrate a smooth homotopy of tangent $n$-planes to construct a regular homotopy immersions of a manifold $N^n$ of dimension $n\ge 2$ \underbar{must fail} since it is well-known that one cannot always integrate, even locally, a smooth family of $n$-planes to obtain an integral $n$-dimensional surface, i.e., whose tangent $n$-planes are in the given family of $n$-planes. There are integrability conditions that must be satisfied in case $n\ge 2$. A totally new idea was required in higher dimensions that would finesse these integrability constraints. This new idea was proposed by Smale.  

Smale's solution was based on a bundle theoretic approach that included an important geometrical ingredient. He proved what is known as the Covering Homotopy Theorem (CHT) for immersions, which for convenience is formulated as follows: Given an immersion of an $n$-disk $f\: D^n\to \bold R^q$, $q\ge n+1$, and a regular homotopy of immersions $F_t\: A\to \bold R^q$, where $A\equiv S^{n-1}\times[0,\ep]\subset D^n$ is a collar neighbourhood of the boundary sphere $S^{n-1}$, such that $F_0=f$ along $A$, then there is a regular homotopy of immersions $G_t\: D^n\to \bold R^q$ such that $G_t=F_t\: B\to \bold R^q$, $t\in[0,1]$, where $B\subset A$ is a smaller collar neighbourhood of $S^{n-1}$, and such that $G_0$ is the given immersion $f\: D^n\to\bold R^q$. Smale's construction of the homotopy extension of immersions $G_t$, $t\in[0,1]$, introduces suitable twistings in the normal directions ($q\ge n+1$) that were explained by Smale (1959) in strictly analytic terms. In order to prove homotopy classification results Smale (1959) requires also that the constructed regular homotopy of immersions $G_t$ be continuous in the data $(f,F_t)$. Smale therefore formulated his CHT as a Serre fibration, equivalent to the following parametric CHT: Let $Z$ be an auxiliary compact manifold. Given a continuous parametrized family of immersions $f(z)\: D^n\to \bold R^q$, $z\in Z$, and a continuous parametrized regular homotopy of immersions $F_t(z)\: A\to \bold R^q$ such that $F_0(z)=f(z)$ along the collar $A$ for all $z\in Z$, then there is a continuous parametrized regular homotopy of immersions $G_t(z)\: D^n\to\bold R^q$, such that for all $(z,t)\in Z\times[0,1]$, $G_0(z)=f(z)$; $G_t(z)=F_t(z)$ along a smaller collar neighbourhood $B\subset A$ of $S^{n-1}$. Employing the parametric CHT and standard bundle theoretic arguments, Smale (1959) proves that the immersions of $S^n$ into $\bold R^q$, $q\ge n+1$, are classified up to regular homotopy by the homotopy group $\pi_n(V_{q,n})$, where $V_{q,n}$ is the Stiefel manifold of all $n$-frames in $\bold R^q$, a space already familiar to homotopy theorists. An interesting consequence of Smale's theory is that the immersions of the 2-sphere $S^2$ into $\bold R^3$ are classified by the group $\pi_2(V_{3,2})=\pi_2(SO(3))=0$. Hence any two immersions of $S^2$ into $\bold R^3$ are regularly homotopic through immersions. This particular result caused a minor sensation in topological circles, and eventually also in the wider mathematical community, because it meant that the standard inclusion of $S^2\subset \bold R^3$ could be connected through a regular homotopy of immersions to the embedding of $S^2$ obtained by reflecting $S^2$ through a 2-plane passing through the centre of the sphere. This was known as an ``eversion'' of the sphere, or more simply, turning the sphere inside out. At that time no explicit geometrical construction of an eversion of $S^2$ was known to Smale or to anyone else. Indeed, over the next few decades several eversions were found which became the subject of separate mathematical papers, an art exhibit that displayed the eversion through a sequence of wire mesh immersed spheres (due to C. Pugh at Berkeley), a movie, and later the subject of an amazing video, Silvio Levy et al.\,(1995), which visually explains Bill Thurston's method of ``corrugations'' for everting the sphere. This video presents a particularly convincing visual proof of the eversion theorem. 

It is essential for Smale's CHT that $q\ge n+1$. Indeed, the CHT does not extend to the equidimensional case $q=n$, i.e., the case of immersions of $D^n$ into $\bold R^n$. Hence the extra dimensions, $q\ge n+1$, provide enough room to prove the CHT, an important insight. To illustrate the problem in the equidimensional case, in the above notation, let $f\: D^2\to\bold R^2$ be the inclusion map and let $F_t\: A\to \bold R^2$ be a regular homotopy of immersions of the collar neighbourhood $A=S^1\times[0,\ep]\subset D^2$ such that $F_0$ is the inclusion and the immersion $F_t$, $t\in[0,1]$, is obtained by a ``finger move'' in $\bold R^2$ which slides a given point $p\in S^1$ along the line joining $p,-p$ so that $F_1(p)=-2p$; outside of a neighbourhood of $p$ the immersion $F_t$, $t\in[0,1]$, is the inclusion map, independent of $t$ (the image $F_1(S^1)$ has two double points near $-p$ and resembles a ``double figure 8''). Elementary arguments prove that in fact the immersion $F_1\: S^1\to \bold R^2$ does not extend to an immersion $G_1\: D^2\to \bold R^2$, so the CHT fails in case $q=n=2$. Similar arguments apply in all cases $q=n\ge 1$. The question of classifying immersions of $S^1$ into $\bold R^2$ that extend to immersions of the disk $D^2$, and in how many inequivalent ways, was solved by S. Blank (cf. Po\'enaru (1968)).        

Smale's theorem on immersions of spheres was quickly generalized by M. Hirsch (1959)
to the classification of immersions of manifolds $N^n$ into $W^q$, $q\ge n+1$, based essentially on Smale's CHT. Attention soon turned to proving classification results for smooth maps $f\: N^n\to W^q$ that satisfied other interesting geometrical and analytical properties. A key problem for each application was the proof of an analogous CHT in each of the new settings. This {\it ad hoc} approach to the CHT was followed throughout the 1960s. Important results here include: the classification of submersions of open manifolds, i.e., smooth maps $f\: N^n\to W^q$, $n\ge q$, of maximal rank $q$, where $N^n$ is {\it open}\,: each connected component of $N^n\setminus\pt N^n$ is non-compact ($\pt N^n$ is the boundary of $N^n$), due to A. Phillips (1967); the classification of maps of rank $\ge k$, known as $k$-mersions, $k< q$, due to S. Feit (1969); the classification of non-degenerate immersed circles in $\bold R^3$, due to E.\,A. Feldman (1968); an analogous immersion theory that applies to classify piecewise-linear immersions of combinatorial manifolds, due to Haefliger and Po\'enaru (1964); Po\'enaru's folding theorem (1966). An important methodological feature throughout this period is the formulation of results in terms of conditions on tangent bundle maps, $Df\: T(N^n)\to T(W^q)$ associated to a smooth map $f\: N^n\to W^q$. This was both a strength and a limitation. A strength because it allowed for easy, geometrical formulations; a limitation because it was not suitable for generalizations to the classification of maps $f\: N^n\to W^q$ that satisfied very general conditions on the $r$th order derivatives of $f$, for all $r\ge 1$. The bundle map formulation is a workable alternative in most cases only when $r=1$. For example, let $\text{Imm}(N^n,W^q)\subset C^{\infty}(N^n,W^q)$ be the open subspace, in the $C^1$ compact-open topology, of immersions of $N^n$ into $W^q$. Let also $\text{Mono}(T(N^n),T(W^q))\subset \text{Hom}(T(N^n),T(W^q))$ be the subspace of tangent bundle homomorphisms consisting of bundle monomorphisms, i.e., which have maximal rank $n$ at each point. Then Hirsch (1959) proves that in case $q\ge n+1$, the derivative map $D\:\text{Imm}(N^n,W^q)\to\text{Mono}(T(N^n),T(W^q))$, $f\mapsto Df$, is a weak homotopy equivalence, where in general a map $g\: X\to Y$ is a {\it weak homotopy equivalence} if the induced map $g_{*}\: \pi_0(X)\to \pi_0(Y)$ is a set bijection on the set of arc components and $g_{*}\:\pi_i(X)\to\pi_i(Y)$ is an isomorphism of homotopy groups for all $i\ge 1$. In particular, (i) at the $\pi_0$ level, the existence of a bundle monomorphism $\alpha$ implies the existence of an immersion $f$ such that $Df,\alpha$ are homotopic through bundle monomorphisms; (ii) at the $\pi_1$ level, if $f,g$ are immersions such that $Df,Dg$ are homotopic through bundle monomorphisms then the immersions $f,g$ are regularly homotopic. In this way the {\it a priori} hard calculation of $\pi_i(\text{Imm}(N^n,W^q))$ is reduced to the easier algebraic calculation of $\pi_i(\text{Mono}(T(N^n),T(W^q))$. Although higher order non-degeneracy conditions on maps $f\: N^n\to W^q$ can be formulated in terms of the symmetric tensor product of a tangent bundle, cf. E.\,A. Feldman (1965), the most general formulation of problems involving arbitrary open and closed conditions on higher order derivatives of $f\: N^n\to W^q$ requires the language of jet spaces, discussed below.

Historically, Thom (1957) provided an influential early exposition in the Bourbaki Seminar of Smale's immersion theory of spheres, based in part on Thom's personal discussions during 1956 with Smale at the University of Chicago, where Thom was then visiting, and where Smale had accepted his first position after completing his thesis under the direction of Raoul Bott at Ann Arbor, Michigan. Thom (1957) proves a general form of the CHT for immersions\,: the restriction of the space of immersions of a manifold to the space of immersions of a submanifold satisfies the covering homotopy property, provided that the immersion of the ambient manifold is in ``good position.'' Thom's constructions, based on good position, provided a more conceptual proof of Smale's CHT above. Thom's good position arguments were developed in courses on immersion theory delivered by Haefliger at Columbia University, and independently by Po\'enaru at Harvard. Po\'enaru's unpublished Harvard course notes of 1964 served also as an inspiration for the theory of submersions of open manifolds developed in Phillips (1967). In particular, Phillips adapts the good position arguments, due originally to Thom in the immersion case, to prove a covering homotopy theorem in the context of submersions of manifolds, a key step in the classification theory of submersions of open manifolds.

Of historical interest also is a refinement of the PL immersion theory of Haefliger and Po\'enaru (1964) made by Lees (1969) in the topological category, where Lees proves a corresponding immersion theorem for topological manifolds. Lees' immersion theorem was then employed by Lashof and Rothenberg (1969) to prove in particular that a compact, simply connected topological manifold $M^m$ of dimension $m\ge 6$, or $m\ge 5$ if $\pt M^m$ is empty, admits a PL manifold structure if $H^4(M^m;\bold Z_2)=0$, which is unique up to PL homeomorphism if $H^3(M^m; \bold Z_2)=0$. Thus topological immersion theory contributed to solving, in part, the classical problem of classifying PL manifold structures on a topological manifold up to PL homeomorphism, completely solved previously only in the low-dimensional cases $m\le 3$. A complete solution to the classification problem for PL structures on topological manifolds in the above dimension range (in particular, no assumption on $\pi_1(M^m)$), finally was obtained by Kirby and Siebenmann (1969), based in particular on Kirby's ``torus trick'' that appeared in Kirby (1969). Only the dimension 4 case remained unsolved.   
\vskip .25cm
\flushpar{\bf 3.\,Jet Spaces.} The following informal presentation of jet spaces is sufficient for our purposes. For a detailed account of jet space bundles in immersion-theoretic topology, cf. Eliashberg and Mishachev (2002), Gromov (1986), Spring (1998). The original formulation of jet spaces is due to C. Ehresmann (1951) with later seminal work in jet spaces due to R. Thom. Given an $r$th order system of P.D.E.s one can replace all the derivatives of order $k$, $1\le k\le r$, of the unknown functions by new independent variables, resulting in systems of equations in these new variables. Thus consider the system of $r$th order P.D.E.s, 
$$
F(x,\pt^{\alpha}f(x))=0\in\bold R^p, \quad |\alpha|\le r, \quad x\in U,\tag1
$$
in the unknown $C^r$-function $f\: U\to \bold R^q$, where $U\subset \bold R^n$ is open and, with respect to coordinates  $(u_1,\dots ,u_n)\in\bold R^n$, the differential operator 
$$\pt^{\alpha}=\frac{\pt^{|\alpha|}}{\pt^{r_1}_{u_1}\circ\pt^{r_2}_{u_2}\circ\dots \circ\pt^{r_n}_{u_n}}\quad\quad\alpha=(r_1,\dots ,r_n),
$$
where $|\alpha| =r_1+r_2+\dots +r_n$; $\pt^0f=f$ if $|\alpha| =0$. Associated to (1) is the system of equations $F(x,p^{\alpha})_{|\alpha|\le r}=0$, where the variables $p^{\alpha}\in\bold R^q$ replace the partial derivatives $\pt^{\alpha}f(x)\in\bold R^q$ in the equation (1). Let $U\subset \bold R^n$, $V\subset \bold R^q$ be open subsets. In terms of the above variables we define the $r$th order jet space 
$$
J^r(U,V)=\{(x,y,p^{\alpha})_{1\le|\alpha|\le r}\mid x\in U, y\in V\}
$$
Thus $J^r(U,V)=U\times V\times\prod\bold R^q$, where the factors $\bold R^q$ correspond bijectively to the differential operators $\pt^{\alpha}$, $1\le |\alpha|\le r$. In particular there is a product bundle $J^r(U,V)\to U\times V$, $(x,y,p^{\alpha})\mapsto (x,y)$, fiber $J^r_{a,b}(U,V)=(a,b)\times\prod\bold R^q$ over $(a,b)\in U\times V$ ($a\in U$ is a ``source'' point; $b\in V$ is a target'' point). Employing Taylor's theorem, each point $w\in J^r(U,V)$ is of the form of an {\it $r$-jet extension}\,: $w =j^rf(x)=(x,f(x),\pt^{\alpha}f(x))_{1\le|\alpha|\le r}$ for some smooth function $f$ defined near $x\in U$, $f(x)\in V$. From this perspective we view $F\: J^r(U,\bold R^q)\to \bold R^p$. A solution to (1) is a $C^r$-function $f\: U\to \bold R^q$ whose $r$-jet extension $j^rf(x)=(x,\pt^{\alpha}f(x))_{|\alpha|\le r}\in F^{-1}(0)\subset J^r(U,\bold R^q)$, for \underbar{all} $x\in U$. 

Gromov (1969) treats in particular open subsets $\Omega\subset J^r_{0,0}(\bold R^n,\bold R^q)$, $r\ge 1$, of interest in topology and geometry. For example $\Omega_k\subset J^1_{0,0}(\bold R^n,\bold R^q)=\prod_1^n\bold R^q$ is the open set defined by the condition that the $q\times n$ matrix, whose columns are the vectors $p^i\in\bold R^q$, has rank $\ge k$, where $p^i$ corresponds to the derivative $\pt/\pt u_i$, $1\le i\le n$, in coordinates $(u_1,\dots ,u_n)\in\bold R^n$. In case $n=k\le q$, $\Omega_n$ defines the immersion relation; in case $n\ge k=q$, $\Omega_q$ defines the submersion relation. A related example, derived from Nash's later work on smooth isometric immersions, is as follows. In terms of the jet space variables $p^i$ above and also the jet space variables $p^{jk}\in\bold R^q$ that correspond to the 2nd partial derivatives $\pt^2/\pt u_j\pt u_k$, $1\le j,k\le n$, the condition that the all the vectors $p^i, p^{jk}\in\bold R^q$, $1\le i,j,k\le n$, are linearly independent defines an open subset, the free map relation $\Omega_{\text{free}}\subset J^2_{0,0}(\bold R^n,\bold R^q)$, $q\ge n+n(n+1)/2$. Because manifolds are locally Euclidean, associated to smooth manifolds $N^n,W^q$ there is an $r$-jet space $J^r(N^n,W^q)$ which is a smooth manifold with charts of the form $J^r(U,V)$, where $U\subset \bold R^n$, $V\subset\bold R^q$ are open. The projection $p^r\: J^r(N^n,W^q)\to N^n\times W^q$, $(x,y,p^{\alpha})_{1\le |\alpha|\le r}\mapsto (x,y)$ is a fiber bundle, fibers of the form $J^r_{0,0}(\bold R^n,\bold R^q)$. Thus associated to an open set $\Omega\subset J^r_{0,0}(\bold R^n,\bold R^q)$ that is invariant under local changes of coordinates there is a subfiber bundle $p^r\: E(\Omega)\to N^n\times W^q$, fiber $\Omega$. In particular the composed projection $p^r\: \Cal R=E(\Omega)\to N^n\times W^q\to N^n$ is a fiber bundle, fiber $W^q\times \Omega$. For example, $\Cal R_k\subset J^1(N^n,W^q)\to N^n$ corresponds to the open set $\Omega_k$ above; in particular $\Cal R_n\subset J^1(N^n,W^q)\to N^n$ is the immersion relation. $\Cal R_{\text{free}}\subset J^2(N^n,W^q)\to N^n$ is the free map relation corresponding to $\Omega_{\text{free}}$. In this context one is interested in $C^r$-maps $f\: N^n\to W^q$ whose $r$-jet extension $j^r f\: N^n\to J^r(N^n,W^q)$, locally of the form, $j^rf(x)=(x,f(x),\pt^{\alpha}f(x))_{1\le |\alpha|\le r}$, satisfies $j^rf(x)\in\Cal R$ for \underbar{all} $x\in N^n$, i.e., the map $f$ {\it solves} the relation $\Cal R$. 

Somewhat more generally, one considers a smooth fiber bundle $p\: X\to N^n$, fiber a manifold $W^q$. Thus locally near each point of the base manifold $N^n$, the map $p$ is smoothly equivalent to a product bundle, fiber $W^q$. An example is the cotangent bundle $T^{*}(N^n)\to N^n$, fiber $\bold R^n$, whose sections are the differential 1-forms on $N^n$. There is an associated $r$-jet space bundle $p^r\: \Xr\to N^n$ whose points are $r$-jet extensions of local sections of the bundle $X\to N^n$. Since local sections of the bundle are maps to the fiber $W^q$, then over a chart $U\subset N^n$, a point $w\in\Xr$ is of the form $w=j^rf(x)$ where $f\in C^r(U,W^q)$, $x\in U$. $\Xr$ is a smooth manifold whose charts are of the form $J^r(U,V)$, where $U\subset \bold R^n$, $V\subset\bold R^q$ are open. In the special case of the product bundle, $X=N^n\times W^q\to N^n$, then $\Xr=J^r(N^n,W^q)$. Associated to an open set $\Omega\subset J^r_{0,0}(\bold R^n,\bold R^q)$ that is invariant under local changes of coordinates there is a subfiber bundle $p^r\: \Cal R=\Cal R(\Omega)\to N^n$, fiber $W^q\times\Omega$.      

Following Gromov (1986), an $r$th order {\it (partial) differential relation} is a subspace $\Cal R$ of an $r$th order jet space $X^{(r)}$. For example, the system (1) of P.D.E.s defines a closed differential relation $\Cal R=F^{-1}(0)\subset J^r(U,\bold R^q)$. A {\it formal solution} to a differential relation $\Cal R$ consists of a continuous section of the jet space bundle $p^r\:\Xr\to N^n$ with values in the differential relation $\Cal R$: a continuous map $g\: N^n\to \Xr$ that is locally of the form, $g(x)=(x,f(x),p^{\alpha}(x))_{1\le|\alpha|\le r}\in \Cal R$, where $f\in C^0(U,W^q)$, $U\subset N^n$ is a chart. In general the jet space variables of a formal solution are far from being the partial derivatives of a function, i.e., we do not assume $p^{\alpha}(x)=\pt^{\alpha}f(x)$ for all $|\alpha|\le r$. In particular, a formal solution  
to the immersion relation consists of a section $g\: N^n\to \Cal R_n\subset J^1(N^n,W^q)$, locally of the form $g(x)=(x,f(x),p^i(x))_{1\le i\le n}$, such that the corresponding vectors $p^i(x)\in \bold R^q$, $1\le i\le n$, are linearly independent. The existence of a formal solution to an open differential relation is a trivial matter locally, while global existence often can be settled by the methods of algebraic topology. The existence of a formal solution to a differential relation is therefore a weak necessary topological condition for solving a differential relation, since in general a formal solution satisfies no integrability conditions on the jet space variables that correspond to higher order partial derivatives.  

A formal solution $g\: N^n\to \Cal R$ is {\it holonomic} if $g$ is a global $r$-jet extension\,: $g=j^rf$ for a $C^r$-section $f\: N^n\to X$\,; locally near each $x\in N^n$, $p^{\alpha}(x)=\pt^{\alpha}f(x)$, for all $|\alpha|\le r$. Thus a holonomic solution is indeed a solution to the given differential relation $\Cal R$. For example, a holonomic solution to the free map relation is a 2-jet extension $j^2f\: N^n\to\Cal R_{\text{free}}\subset J^2(N^n,W^q)$, for some $C^2$-map $f\: N^n\to W^q$. Thus locally at each point $x\in N^n$, the derivatives $\tfrac{\pt f}{\pt x_i}(x),\tfrac{\pt^2 f}{\pt x_j\pt x_k}(x)$, $1\le i,j,k\le n$, are linearly independent vectors in $\bold R^q$. Following Gromov, a differential relation satisfies the {\it $h$-principle} (homotopy principle, coined by Gromov (1986)) if every formal solution is homotopic through formal solutions to a holonomic solution. There is also a parametric $h$-principle in case of parametrized families of formal solutions. All of the classical problems and results in immersion-theoretic topology developed by Smale and later researchers can be reformulated in terms of solutions to the $h$-principle and also to the parametric $h$-principle in various jet bundle contexts. The surprising discovery of all this research is that the existence of a formal solution, a seemingly weak topological condition involving no integrability conditions, is in fact sufficient for the existence of holonomic solutions in the very broad contexts in which the $h$-principle was proved to hold. This unexpected insight of immersion-theoretic topology constitutes one of its major achievements. 

Around the time of Smale's classification of immersions of spheres in Euclidean space, it was clear to Ren\'e Thom in France that Smale's work should be recast in terms of jet spaces. Indeed, Thom had worked on the theory of singularities of maps for several years, in which the language of jet spaces was natural. Thom (1959) was the first to propose a general jet space formulation of the problems to be studied in immersion-theoretic topology. Interestingly, Thom states his classification results in terms of certain homology equivalences rather than in terms of weak homotopy equivalences that were employed in the immersion case by Smale (1959) and in other special cases discussed above by later researchers. However, Thom was unable to prove significant results (cf. Smale's review in MR {\bf 22}\,\#12537) because, with hindsight, he did not have available to him the relevant geometrical tools. These geometrical and analytical tools, in the jet space context, would begin to appear a decade later in the work of Gromov and Eliashberg at Leningrad. In general however, the jet space formulation of problems in immersion-theoretic topology was not employed by topologists in the U.S. who were working in this area. To these topologists the jet space formulation seemed unnecessary and it obscured the tangent bundle formulation which seemed more natural. Thus it was primarily the European school of topologists, including Thom, Po\'enaru, Haefliger, who were receptive to the importance of jet spaces for studying general problems in immersion-theoretic topology. This would prove to be important for the developments that would soon take place in Rokhlin's seminar in Leningrad. 
\vskip .25cm
\flushpar{\bf 4.\,The Leningrad School.} As a graduate student under the direction of V. Rokhlin, Misha Gromov undertook to reexamine the foundations of Smale's immersion theory and the related papers by Phillips (1967) and others on this topic, in order to formulate an independent view of the whole subject. In a first paper, Gromov (1968) considers smooth maps $f$ of a manifold $M^m$, equipped with a foliation $\Cal F$ (cf.\,\S5 for foliations), to a manifold $N^n$ equipped with a subbundle $\xi$, such that on each leaf $\Cal L$ of $\Cal F$, the tangent map $df\: T(\Cal L)\to T(N^n)/\xi$ is injective on each fiber. In case $\dim\Cal F< \text{codim\,}\xi=n-\dim\xi$, Gromov proves an appropriate weak homotopy equivalence result, which recovers as a special case the Smale-Hirsch theory of immersions of manifolds, for which $\Cal F=T(M)$, $\xi =0$, and $n\ge m+1$. There soon followed Gromov's thesis (1969), the striking feature of which is Gromov's reformulation of the entire subject in the language of sheaves, with applications to the context of jet spaces. As explained above, this application to jet spaces was closer to the point of view of the French school, especially the work of R. Thom, rather than the U.S. school of topology. Gromov told me that he had learned about jet spaces from R. Palais (1965, ch.\,IV), a book about the Atiyah-Singer index theorem. The key advance made by Gromov's thesis (1969) was to prove a \underbar{universal} version of Smale's CHT in a very general sheaf-theoretic setting that applied, in a special case, to the context of suitable open differential relations $\Cal R\subset X^{(r)}$ in jet space bundles such that $p^r\: \Cal R\to N^n$ is a subfiber bundle, where $N^n$ is an open manifold. In particular the theory in Gromov (1969) applies to subbundles $p^r\:\Cal R=\Cal R(\Omega)\to N^n$ associated to an open subset $\Omega \subset J^r_{0,0}(\bold R^n,\bold R^q)$ as discussed above. For suitable choices of $\Cal \Omega$ one obtains the immersion relation, the submersion relation, the exact symplectic form relation, and many other classical relations. Gromov also requires an additional ``naturality'' condition, satisfied in all examples of interest, that allows one to lift diffeomorphisms of the base manifold $N^n$ to fiber preserving diffeomorphisms of the fiber bundle $X\to N^n$, and hence to fiber preserving diffeomorphisms of the bundle $p^r\:\Xr\to N^n$. The relation $\Cal R\subset\Xr$ is {\it natural} if these lifted diffeomorphisms leave $\Cal R$ invariant. These lifts are employed primarily to ensure that the required twisting constructions in Gromov's universal CHT, analogous to those used by Smale in the immersion case, leave invariant these open differential relations $\Cal R$. With these innovations, he followed the Smale program in terms of towers of bundles to prove homotopy classification theorems. In somewhat more technical terms, Gromov proved the following. Let $i\: \Cal H\to\Cal F$ be the inclusion map of the space of holonomic solutions into the space of formal solutions, in the compact-open topology, with respect to a given partial differential relation $\Cal R$. If $p^r\:\Cal R\to N^n$ is a subfiber bundle where the base manifold $N^n$ is open and the differential relation $\Cal R\subset\Xr$ is both open and satisfies the above naturality condition then Gromov (1969) proves that the inclusion map $i$ is a weak homotopy equivalence. The $h$-principle follows from the bijection $i_{*}\: \pi_0(\Cal H)\to\pi_0(\Cal F)$. The parametrized $h$-principle follows from the vanishing of the relative homotopy groups $\pi_i(\Cal F,\Cal H)=0$, for all $i\ge 1$. Gromov (1969) recovered as special cases most of the results that had been obtained to date and much more, including existence on open manifolds of exact symplectic forms, contact structures, foliations, and the somewhat surprising result that on an open manifold there exists a Riemannian metric (not complete in general) for which the sectional curvature is always $>0$ (or always $<0$). Gromov's thesis (1969) provides an explicit analytic algorithm for proving his universal sheaf-theoretic CHT,  based on his beautiful and surprisingly simple picture of the whole process of the CHT. Although his picture was not published it has been communicated informally in topological circles ever since (3 decades!). Gromov's universal CHT incorporates the relevant geometrical features of the CHT proved in the immersion case by Smale (1959) and in the submersion case by Phillips (1967). Gromov's thesis was brought to the West by Tony Phillips who met with Gromov in Leningrad around 1970. The essentials of Gromov's thesis were published in the West as a set of expository notes by Haefliger (1971), and also Po\'enaru (1971) who provided some pictorial details. An interesting refinement to Gromov's CHT was proved by Phillips (1974) who classifies maps $f\: V^n\to W^q$ of constant rank $k$, where $V^n$ is an open manifold that admits a proper Morse function with no critical points of index $> k$. This result does not follow from Gromov (1969) since the differential relation that corresponds to maps of constant rank is a \underbar{closed} condition, not an open condition as required in Gromov (1969). For this purpose Phillips proves a ``weak micro-covering homotopy lemma'' that generalizes Gromov's CHT in the special case of maps of constant rank (independently, using different methods, Gromov (1973a) proved a similar result for maps of constant rank).       

We comment here on the historical evolution towards employing open base manifolds $N^n$ in Gromov (1969). Note that an immersion of $S^1$ into $\bold R^2$ extends in the normal directions to an immersion ($=$ submersion) of an annulus $S^1\times(-\ep,\ep)$ into $\bold R^2$. Conversely, a submersion of the annulus $S^1\times(-1,1)$ into $\bold R^2$ restrict to an immersion of $S^1=S^1\times \{0\}$ into $\bold R^2$. More generally, it was pointed out in Phillips (1967) that the Hirsch (1959) classification of immersions of $N^n$ into $W^q$, $q\ge n+1$, can be recovered from Phillips (1967) applied to the submersions of an open $q$-dimensional manifold $U^q$ into $W^q$, where $U^q$ is a neighbourhood of $N^n$ in a suitable normal $(q-n)$-bundle over $N^n$. Secondly, employing elementary topological arguments, an open manifold $N^n$ is diffeomorphic to a small neighbourhood of a subcomplex $K$ of $N^n$, $\dim K\le n-1$ (called a spine of the open manifold $N^n$), in a suitable triangulation of $N^n$. In particular, each simplex of $K$ has a normal bundle of dimension $\ge 1$ in $N^n$. Thus if $N^n=S^1\times (-1,1)$ let $K=S^1\times\{0\}$. Po\'enaru (1962), Hirsch (1961) took advantage of these normal directions to simplices in $K$ to prove some refinements to immersion theory in the case of open manifolds. For example these authors prove that an open manifold $N^n$ can be immersed in $\bold R^n$ if $T(N^n)$ is a trivial $n$-plane bundle. Phillips (1967) uses a similar idea, but he replaces triangulations with suitable handlebody decompositions of an open manifold $N^n$ that have handles only of dimension $\le n-1$. Thus if $N^n$ is open it is sufficient to prove the appropriate CHT in a small neighbourhood of $K$. In particular, the suitable twistings required by the CHT can take place in the normal directions to these simplices of dimension $\le n-1$ in $K$, and hence are \underbar{internal} constructions in $N^n$, not ambient constructions in $W^q$. These ideas were incorporated in Gromov (1969) whose proof of a universal CHT employs suitable twistings within the open base manifold $N^n$ in normal directions to simplices in a spine $K$. By contrast, if $N^n$ is {\it closed} (compact without boundary), then a lower-dimensional spine $K\subset N^n$ does not exist. Hence the CHT will apply over simplices (or handles) of dimension $\le n-1$, but, for lack of normal directions, will generally not apply over the top-dimensional $n$-simplices (or $n$-handles) of $N^n$, which in turn traces back to the failure of Smale's CHT in the equidimensional case $q=n$ (cf.\,\S2). Thus open manifolds are natural domains for covering homotopy proofs of the $h$-principle. Nevertheless if $N^n$ is closed then the $h$-principle is provable for open relations $\Cal R$ that admit ``microextensions'', i.e., that extend locally to suitable open relations over $N^n\times \bold R^m$, for some $m\ge1$, (thereby obtaining extra normal directions to apply an adapted CHT). Cf.\,Gromov (1986 \S2.2.4) for the general theory with applications (including Feit's $k$-mersion theory for closed manifolds), introduced as ``micro-majorization'' in Gromov (1972, \S6.1.7; \S6.2.6); also du\,Plessis (1975).  

The Leningrad school of immersion-theoretic topology gave rise to three general methods for solving differential relations in jet spaces: (i) the covering homotopy method discussed above in Gromov's thesis, also known as the method of continuous sheaves; (ii) the method of removal of singularities due to Gromov and Eliashberg (1971a), which was their first collaborative work; (iii) the method of convex integration due to Gromov (1973b). These general methods are not linearly related in the sense that successive methods subsumed the previous methods. Each method had its own distinct foundation, based on an independent geometrical or analytical insight. Consequently each method has a range of applications to problems in topology that are best suited to its particular insight. For example, the method of removal of singularities uniquely applies also to the complex analytic and algebraic setting, cf. Gromov and Eliashberg (1971b), (1992). Convex integration theory, discussed below, applies to solve closed differential relations in jet spaces, including certain general classes of underdetermined non-linear systems of partial differential equations. On the other hand many classical results in immersion-theoretic topology, such as the classifications of immersions, are provable by all three methods. The book Adachi (1993) provides an exposition of Gromov (1969) and Gromov (1973b).             

Another important contribution to immersion-theoretic topology from the Leningrad school at that time was the thesis of Yasha Eliashberg (1972) on the surgery of singularities of maps. Eliashberg was well-acquainted with the main results and methods of immersion-theoretic topology from the West and also with Gromov's recent innovations. Eliashberg turned his attention to simplifying the {\it singularities} (the locus of points where the rank is not maximal) of maps $f\: N^n\to W^q$, $n\ge q$, $N^n$ compact, as much as homotopy theory would allow. For example a smooth map $f\: N^n\to\bold R^q$, $n\ge q$, where $N^n$ is compact without boundary, must have singularities. Indeed, if $f$ had maximal rank $(=q)$ at each point then, by the implicit function theorem, the image of $f$ is an open set in $\bold R^q$ which is impossible since $N^n$ is compact.  Thus in cases where singularities cannot be avoided, the problem was to determine the simplest forms of these singularities, up to a pointwise small perturbation. For this task, the jet space formulation was essential and the results of the French school of singularity theory, which in turn was inspired by the fundamental results of H. Whitney (1955) on generic singularities of mappings into the plane, played an important role. Eliashberg (1970, in the case $n=q$), (1972, in the case $n\ge q\ge 2$) proved that, subject to mild bundle hypotheses, up to pointwise small perturbations, one could arrange that maps only had the simplest of singularities: fold singularities that occured along prescribed submanifolds in $N^n$. Explicitly, a smooth map $f\: N^n\to W^q$, $n\ge q$, has a fold singularity of index $s$ at $x\in N^n$ if near $x$ the map $f$ is equivalent to the map, 
$$
\pi\: \bold R^n =\bold R^{q-1}\times \bold R^{n-q+1}\to\bold R^q =\bold R^{q-1}\times \bold R, 
$$
given by the formula
$$
(y,x_1,\dots ,x_{n-q+1})\mapsto \left(y,-\sum_1^sx_i^2+\sum_{s+1}^{n-q+1}x_j^2\right),
$$
where $y\in\bold R^{q-1}$. In case $q=1$ the $y$-term disappears and we recover the classical quadratic Morse singularities of index $s$ ot the origin in $\bold R^n$. In this model the singularities of $\pi$ occur along the submanifold $\bold R^{q-1}\equiv\bold R^{q-1}\times\{0\}\subset \bold R^n$, at which points the map $\pi$ has rank $q-1$ and has quadratic Morse singularities in the normal directions $\{y\}\times\bold R^{n-q+1}$, $y\in\bold R^{q-1}$. Subject to mild bundle hypotheses Eliashberg (1972) proved that up to a $C^0$-small perturbation, i.e., a pointwise small perturbation, one may assume that the map $f$ is a submersion (maximal rank $=q$) on $N^n\setminus\bigcup_1^k V_i$, where $V_i$ is a prescribed $(q-1)$-submanifold of $N^n$, $q\ge 2$, and such that along each $V_i$ the map $f$ locally has the form of a fold singularity as above, $1\le i\le k$. Previous results of this type, in particular eliminating cusp singularities, were obtained by Levine (1965) in the special case $q=2$. These results on the topological simplification of singularities were in marked contrast to the results on generic singularities obtained by Thom, and by the Arnold school of singularity theory in Moscow. Thus Arnold was interested in determining the form of generic singularities in different contexts, up to $C^{\infty}$-small perturbations, which led to his intricate local classifications of singularities. Eliashberg showed that if you relax the perturbations to be only $C^0$-small, then fold singularities are sufficient for topological purposes. These results formed the basis of his more recent work, in collaboration with Kolya Mishachev, that simplifies singularities of parametrized families of maps to have only fold and cusp singularities of prescribed type, called wrinkled singularities, cf. Eliashberg and Mishachev (1997), (2000); also Spring (2002).

We conclude our survey of the contributions to immersion-theoretic topology from the Leningrad school with the article Gromov (1973b) on convex integration theory. Gromov's thesis (1969) had clarified, in a general sheaf-theoretic formulation, the Smale CHT and its role in immersion-theoretic topology. Gromov now began thinking about how to exploit the geometrical and algebraic structure of jet spaces themselves, from a topological point of view, in order to formulate constructions in a jet space context that would lead to a proof of the $h$-principle. In this respect he was moving into uncharted territory and he was way ahead of his times. Although his 1973 paper dealt only with the $h$-principle in the case of differential relations in the space of 1-jets it was clear that Gromov had further applications in mind in the general context of spaces of $r$-jets for all $r\ge 1$. There were several striking features of his 1973 paper. First, he introduced a new analytical approximation lemma (the one-dimensional lemma) which replaced the use of Smale's CHT. Thus Gromov had discovered a new proof procedure for proving the $h$-principle that circumvented both Smale's CHT and also the tower of bundles approach initiated by Smale that was based on the CHT. Furthermore this proof procedure worked whether or not the base manifold was compact or was an open manifold, the latter being a necessary hypothesis in his thesis. Second, the method of convex integration applied for the first time to solve certain classes of \underbar{closed} differential relations in jet spaces. In particular, Gromov introduced a convergence procedure that solved certain classes of non-linear systems of P.D.E.s that satisfied some mild local convexity properties. The solution of the $h$-principle for general closed differential relations in jet spaces had long remained an unsolved puzzle. Hence Gromov (1973b) was a major advance in immersion-theoretic topology with respect to the solution of closed differential relations. 

The basic content of Gromov's one-dimensional lemma is as follows. Let $A\subset\bold R^q$ be open and connected, and let $f\in C^1([0,1],\bold R^q)$, $q\ge 2$ such that for all $t\in[0,1]$, $f^\prime(t)\in\text{Conv}(A)$, the convex hull of $A$ in $\bold R^q$. Let also $\ep>0$. There exists $g\in C^1([0,1],\bold R^q)$ such that for all $t\in[0,1]$: (i) $|g(t)-f(t)|\le \ep$; (ii) $g^\prime(t)\in A$. Thus the map $f$ whose derivative $f^\prime$ has values in $\text{Conv}(A)$ is replaced by a $C^0$-approximation $g$ whose derivatives $g^\prime$ takes its values in $A$. In general $A$ could be a thin open set with a large convex hull. For applications of immersion-theoretic topology to $r$th order jet spaces there is also a $C^r$-version of the one-dimensional lemma for all $r\ge 1$. Let $f\in C^r(I^n,\bold R^q)$, $I^n=[0,1]^n$, such that, with respect to coordinates $x=(u_1,\dots ,u_{n-1},u_n=t)\in \bold R^n$, $\tfrac{\pt^r f}{\pt t^r}(x)\in\text{Conv}(A)$ for all $x\in I^n$. Let $\ep>0$. There exists $g\in C^r(I^n,\bold R^q)$ such that for all $x\in I^n$: (i) $\tfrac{\pt^rg}{\pt t^r}(x)\in A$; (ii) $|\pt^{\alpha}(f-g)(x)|\le \ep$ for all derivatives $\pt^{\alpha}$ such that $|\alpha|\le r$ and $\pt^{\alpha}\ne \pt^r/\pt t^r$. In particular $g$ is a $C^{r-1}$-small approximation to $f$. These $C^{r-1}$-approximation results for $r\ge 2$ were a vast improvement over the Smale theory, as developed in Gromov's thesis (1969), which was only a $C^0$-approximation theory in general, even in the context of differential relations in $r$th order jets spaces, $r>1$. To illustrate, let $f\: N^n\to W^q$ be an immersion whose 1-jet extension $j^1f\: N^n\to J^1(N^n,W^q)$ can itself be extended to a formal solution of the free map relation.  Thus there is a formal solution $g\: N^n\to \Cal R_{\text{free}}\subset J^2(N^n,W^q)$ such that locally at each $x=(x_1,\dots ,x_n)\in N^n$,
$$
g(x)=(j^1f(x),p^{jk}(x))=(x,f(x),\pt_{x_i}f(x),p^{jk}(x))\in\Cal R_{\text{free}}\quad 1\le i,j,k\le n,
$$
where the vectors $\pt_{x_i}f(x), p^{jk}(x)\in\bold R^q$, $1\le i,j,k\le n$ are linearly independent. In the extra dimensional case $q> n+n(n+1)/2$, convex integration theory proves a $C^1$-dense $h$-principle: the formal solution $g$ is homotopic through formal solutions of $\Cal R_{\text{free}}$ to a holonomic solution $j^2 h\: N^n\to\Cal R_{\text{free}}$ such that in addition the immersions $f,h$ are $C^1$-close. This $C^1$-approximation result is beyond the scope of the CHT method of Gromov (1969), but does overlap with results proved earlier by the method of removal of singularities in Gromov and Eliashberg (1971a). Strictly speaking Gromov's approximation lemma, in the context of $C^r$-functions for all $r\ge 1$, was published much later in his tome on the whole of immersion-theoretic topology up to that point, Gromov (1986). In addition Gromov (1973b) states that he can reprove the Nash $C^1$-isometric immersion theorem by the methods of convex integration. Indeed, according to Gromov, the Nash twisting process in normal directions (which serves locally to stretch successive first order partial derivatives) is a kind of one-dimensional lemma. In this sense the Nash paper (1954) is an important precursor to convex integration theory. 

The results of convex integration theory included many of the results of his thesis, and also many new results, especially results on the $h$-principle applied to embeddings\,: the perturbations to yield a holonomic solution had to proceed by  isotopies of embeddings rather than through regular homotopies of immersions. These latter results were explicated only much later in his book, Gromov (1986, \S2.4.5(C)); also Spring (1998, \S8.4.4), Eliashberg and Mishachev (2002, \S4.4-\S4.6). 

\flushpar{\bf 5.\,Foliations.} One of the more interesting applications of immersion-theoretic topology during this period was the classification of foliations on open manifolds. A smooth foliation of dimension $p$ on a manifold $N^n$, $0\le p\le n$, consists of a decomposition of $N^n$ into a disjoint union of connected subsets $(L_{\alpha})_{\alpha\in A}$ (known as the ``leaves'' of the foliation) such that in local coordinates $(x_1,x_2,\dots ,x_n)$ in a suitable smooth chart $U$ at each point of $N^n$ the intersection of the leaves with $U$ is determined by the ``slice'' condition: $x_j=c_j$, $p+1\le j\le n$, where $(c_{p+1}, \dots ,c_n)\in\bold R^{n-p}$. Thus in these coordinates a leaf $L_{\alpha}$ meets $U$ in a possibly infinite number of parallel $p$-dimensional affine subspaces (cf. Lawson (1974) for a fuller discussion with examples). Thus a product manifold $P^p\times Q^q$ determines a foliation of dimension $p$ whose leaves are the slices $P^p\times \{x\}$, $x\in Q^q$. Note that the family of tangent $p$-planes to all of the leaves $L_{\alpha}$ defines an integrable $p$-dimensional subbundle of $T(N^n)$. The basic problem studied was the following form of the $h$-principle. Let $E\subset T(N^n)$ be a $p$-dimensional subbundle: a family of tangent $p$-planes $\{E_x\subset T_x(N^n)\mid x\in N^n\}$ that depends continuously on $x\in N^n$. The $h$-principle requires the existence of a homotopy of $p$-dimensional subbundles of $T(N^n)$ that connects $E$ to an integrable $p$-dimensional subbundle, i.e., to a foliation of $N^n$. In case $N^n$ is an open manifold the immersion-theoretic techniques introduced by Gromov (1969), Phillips (1969), (1970), solved this $h$-principle in important special cases, with a comprehensive formulation provided by Haefliger (1970). Employing more refined geometrical techniques, the $h$-principle for foliations on closed manifolds was solved by Thurston (1974) for codimensions $>1$ and Thurston (1976) for codimension 1 foliations. In particular, $(n-1)$-dimensional foliations exist on a closed $n$-dimensional manifold if the Euler characteristic is zero. Thurston is a co-recipient of the Fields medal in 1982 for which his work above on foliations is included in the citation.             

\flushpar{\bf 6.\,Postscript.} Throughout the history of immersion-theoretical topology there has been a constant struggle to clarify the geometrical constructions that enable one to prove the $h$-principle in increasingly general contexts, with the widest possible scope of applications to topology and geometry. The original Smale CHT was clarified in Gromov's thesis (1969). The general exposition of immersion-theoretic topology, Gromov (1986, Ch.\,2.2, 2.4), further refines and generalizes the covering homotopy methods of Gromov (1969) and also the convex integration theory in Gromov (1973b). More recently Eliashberg and Mishachev (2001), (2002) have proved a holonomic approximation theorem that captures the main ingredients of some previous geometrical constructions, and which appears to put the main geometrical ideas initiated by Smale (1959), later generalized by Gromov (1969), and also Gromov (1986, Ch.\,2.2), into a definitive form. In other directions, Bierstone (1974) proved an equivariant version of the main jet space results in Gromov's thesis (1969), with recent improvements by Datta and Mukherjee (1998). Ralph Cohen (1985) proved the long-standing conjecture that every $n$-dimensional manifold can be immersed in $\bold R^{2n-\alpha(n)}$, where $\alpha(n)$ is the number of ones in the binary expansion of $n$ (this upper bound is best possible). In the case of manifolds that have additional structure, e.g. symplectic or contact structure, then the $h$-principle can be proved for immersions and embeddings that preserve that structure; cf.\,Gromov (1986, \S3.4), Eliashberg and Mishachev (2002, Ch.\,12). From this perspective, one can appreciate the long and fruitful history of immersion-theoretic topology initiated by Steve Smale over 40 years ago. 
\vskip .25cm
\flushpar{\bf References}
\vskip .25cm
\eightpoint{
\item{}Adachi, M., (1993) Embeddings and immersions. Transl. Math. Monographs Vol. 124, Amer. Math. 
\itemitem{}Soc., Providence R.I. MR {\bf 95a:}\,57039.
\item{} Bierstone, E., (1974) Equivariant Gromov Theory. Topology {\bf 13}, 327-345. MR {\bf 50}\,\#5818.
\item{}Cohen, R.L., (1985) The immersion conjecture for differentiable manifolds. Ann.\,of Math. {\bf 122}, 237-328. 
\itemitem{}MR {\bf 86m:}\,57030.
\item{} Datta, M., Mukherjee, A., (1998) Parametric homotopy principle of some partial differential relations.  
\itemitem{}Math. Slovaca {\bf 48}, 411-421. MR {\bf 2000f:}\,58021.
\item{}Ehresmann, C., (1951) Les prolongements d'une vari\'et\'e diff\'erentiable.\,I, II, III. C.R.\,Acad.\,Sci.\,Paris. \itemitem{}{\bf 233}, 598-600; 777-779; 1081-1083. MR {\bf 13},\,386a; {\bf 13},\,584b; {\bf 13},\,584c. Reprinted in\,:\,Charles Ehresmann\,: \oe uvres compl\`etes et comment\'ees. I-1, I-2. Topologie alg\'ebrique et g\'eom\'etrie diff\'erentielle. Cahiers Topologie et G\'eom.\,Diff\'erentielle.\,{\bf 24} (1983), suppl.1 (1984), 343-345; 346-348; 349-351. MR {\bf 86i:}\,01059.
\item{}Eliashberg, Ja.M., (1970) Singularities of folding type. Izv. Akad. Nauk SSSR Ser. Mat. {\bf 34}, 
\itemitem{}1110-1126 (Russian). English translation: Math. USSR Izvestia {\bf 4} (1970), 1119-1134. MR {\bf 43}\,\#4051. 
\item{}Eliashberg, Ja.M., (1972) Surgery of singularities of smooth mappings. Izv. Akad. Nauk SSSR Ser. 
\itemitem{}Mat. {\bf 36}, 1321-1347 (Russian). English translation: Math. USSR Izvestia {\bf 6} (1972), 1302-1326. 
MR {\bf 49}\,\#4021.
\item{}Eliashberg, Y., Mishachev, N., (1997) Wrinkling of smooth mappings and its applications.\,I. Invent.  
\itemitem{}math. {\bf 130}, 345-369. MR {\bf 99d:}\,57021.
\item{}Eliashberg, Y., Mishachev, N., (2000) Wrinkling of smooth mappings.\,II. Wrinkling of embeddings and  
\itemitem{}K. Igusa's theorem, Topology {\bf 39}, 711-732. MR {\bf 2001g:}\,57058.
\item{}Eliashberg, Y., Mishachev, N., (2001) Holonomic approximation and Gromov's $h$-principle. Essays on 
\itemitem{} geometry and related topics. Monogr. Enseign. Math. {\bf 38}, Vol.\,1,2. 271-285. 
\item{}Eliashberg, Y., Mishachev, N., (2002) Introduction to the $h$-principle. Grad. Stud. in Math. Vol. 48, 
\itemitem{}Amer. Math. Soc., Providence, R.I.
\item{}Feldman, E.A., (1965) The geometry of immersions.\,I. Trans. Amer. Math. Soc. {\bf 120}, 185-224. 
\itemitem{}MR {\bf 32}\,\#3065.
\item{}Feldman, E.A. (1968) Deformations of closed space curves. J. Diff. Geom. {\bf 2}, 67-75. MR {\bf 38}\,\#725. 
\item{}Feit, S.D., (1969) $k$-mersions of manifolds. Acta Math. {\bf 122}, 173-195. MR {\bf 39}\,\#4862.
\item{}Gromov, M.L., (1968) Transversal mappings of foliations. Dokl. Akad. Nauk SSSR. {\bf 182}, 255-258 
\itemitem{}(Russian). English translation: Soviet Math. Dokl. {\bf 9} (1968), 1126-1129. MR {\bf 38}\,\#6628. 
\item{}Gromov, M.L., (1969) Stable maps of foliations into manifolds. Izv. Akad. Nauk SSSR Ser. Mat. {\bf 33},
\itemitem{}707-734 (Russian). English translation: Trans. Math. USSR Izvestia {\bf 3} (1969), 671-694. MR {\bf 41}\,\#7708.
\item{}Gromov, M.L., (1972) Smoothing and inversion of differential operators. Math.\,Sbornik {\bf 88(130)},
\itemitem{}382-441 (Russian). English translation: Math.\,USSR Sbornik {\bf 17} (1972), 381-435. MR {\bf46}\,\#10022.
\item{}Gromov, M.L., (1973a) Degenerate smooth mappings. Mat. Zametki {\bf 14}, 509-516 (Russian). English 
\itemitem{}translation: Math. Notes {\bf 14} (1973) 849-853. MR {\bf 49}\,\#3972.
\item{}Gromov, M.L., (1973b) Convex integration of differential relations.\,I. Izv. Akad. Nauk SSSR Ser. Mat.  
\itemitem{}{\bf 37}, 329-343 (Russian). English translation: Math.USSR Izvestia {\bf 7} (1973), 329-343. MR {\bf 54}\,\#1323.
\item{}Gromov, M.L., (1986) Partial differential relations. Ergebnisse der Mathematik und ihrer Grenzgebiete 
\itemitem{}3. Folge $\cdot$ Band 9. Springer-Verlag, Berlin. MR {\bf 90a:}\,58201.
\item{}Gromov, M.L., Eliashberg, Ja.M., (1971a) Removal of singularities of smooth mappings. Izv. Akad. 
\itemitem{}Nauk SSSR Ser. Mat. {\bf 35}, 600-626 (Russian). English translation: Math. USSR Izvestia {\bf 5} (1971), 615-639. MR {\bf 46}\,\#903.
\item{}Gromov, M.L., Eliashberg, Ja.M., (1971b) Nonsingular maps of Stein manifolds. Funkcional.\,Anal.\,i 
\itemitem{}Prilozhen, {\bf 2}, 82-83 (Russian). English translation: Functional Anal. Appl. {\bf 5} (1971), 156-157. MR {\bf 46}\,\#394.
\item{}Gromov, M., Eliashberg, Y., (1992) Embedding of Stein manifolds of dimension $n$ into the affine  
\itemitem{}space of dimension $(3n/2)+1$. Ann.\,of Math.\,(2) {\bf 136}, 123-135. MR {\bf 93g:}\,32037.
\item{}Haefliger, A., Po\'enaru, V., (1964) La classification des immersions combinatoires. Inst.\,Hautes \'Etudes 
\itemitem{}Sci.\, Publ. Math. {\bf 23}, 75-91. MR {\bf 30}\,\#2515.
\item{}Haefliger, A., (1970) Feuilletages sur les vari\'et\'es ouvertes. Topology {\bf 9}, 183-194. MR {\bf 41}\,\#7709.
\item{}Haefliger, A., (1971) Lectures on the theorem of Gromov. Proceedings of Liverpool Singularities Sym- 
\itemitem{}posium,\,II (1969/1970).\,pp.\,128-141.\,Lecture Notes in Math.\,Vol 209.\,Springer. MR {\bf 48}\,\#12560. 
\item{}Hirsch, M.W., (1959) Immersions of manifolds. Trans. Amer. Math. Soc. {\bf 93}, 242-276. MR {\bf 22}\,\#9980.
\item{}Hirsch, M.W., (1961) On imbedding differentiable manifolds in euclidean space. Ann.\,of Math. {\bf 73}, 
\itemitem{}566-571. MR {\bf 23}\,\#A2223.
\item{} Hirsch, M.W., (1993) The work of Stephen Smale in differential topology. From Topology to Computa-
\itemitem{}tion: Proceedings
of the Smalefest (M.W. Hirsch, J.E. Marsden, M. Shub, eds.) Springer Verlag, New York, 83-106. MR {\bf 97b:}\,57001.
\item{}Kirby, R.C., (1969) Stable homeomorphisms and the annulus conjecture. Ann.\,of Math. {\bf 89} 575-582.
\itemitem{}MR {\bf 39}\,\#3499.
\item{}Kirby, R.C., Siebenmann, L.C., (1969) On the triangulation of manifolds and the Hauptvermutung. 
\itemitem{}Bull.\,Amer.\,Math.\,Soc. {\bf 75}, 742-749. MR {\bf 39}\,\#3500.
\item{}Kuiper, N.H., (1955) On $C^1$-isometric imbeddings.\,I, II. Indag. Math. {\bf 17}, 545-556, 683-689. 
\itemitem{}MR {\bf 17},\,782c.
\item{}Lashof, R., Rothenberg, M., (1969) Triangulation of manifolds.\,I, II. Bull.\,Amer.\,Math.\,Soc. {\bf 75}, 
\itemitem{}750-754; 755-757. MR {\bf 40}\,\#895.
\item{}Lawson, H.B., (1974) Foliations. Bull.\,Amer.\,Math.\,Soc. {\bf 80}, 369-418. MR {\bf 49}\,\#8031.
\item{}Lees, J.A., (1969) Immersions and surgeries of topological manifolds. Bull.\,Amer.\,Math.\,Soc. {\bf 75}, 
\itemitem{}529-534. MR {\bf 39}\,\#959.
\item{}Levy, S., Maxwell, D., Munzner, T., (1995) Outside In (Video). Accompanying text, {\it Making waves},
\itemitem{}University of Minnesota Geometry Center, Minneapolis MN. Distributed by A.K. Peters Ltd., Wellesley, MA. 02181 U.S.A. MR {\bf 96m:}\,00012.
\item{}Levine, H.I., (1965) Elimination of cusps. Topology {\bf 3}, suppl.\,2, 263-296. MR {\bf 31}\,\#756.
\item{}Nash, J., (1954) $C^1$-isometric imbeddings. Ann.\,of Math. {\bf 63}, 383-396. MR {\bf 16},\,515e.
\item{}Palais, R.S., (1965) Seminar on the Atiyah-Singer index theorem. Annals of Math.\,Studies. Vol.\,57.  
\itemitem{}Princeton University Press. MR {\bf 33}\,\#6649.
\item{}Phillips, A., (1967) Submersions of open manifolds. Topology {\bf 6}, 171-206. MR {\bf 34}\,\#8420.
\item{}Phillips, A., (1969) Foliations on open manifolds.\,II. Comment.\,Math.\,Helv. {\bf 44}, 367-370. MR {\bf 40}\,\#6579.
\item{}Phillips, A., (1970) Smooth maps transverse to a foliation. Bull.\,Amer.\,Math.\,Soc., {\bf 76}, 792-797. 
\itemitem{}MR {\bf 41}\,\#7711.
\item{}Phillips, A., (1974) Smooth maps of constant rank. Bull.\,Amer.\,Math.\,Soc. {\bf 80}, 513-517. 
\itemitem{}MR {\bf 48}\,\#12557.
\item{}du\,Plessis, A., (1975) Maps without certain singularities. Comment.\,Math.\,Helv. {\bf 50}, 363-382.
\itemitem{} MR {\bf 53}\,\#1637.
\item{}Po\'enaru, V., (1962) Sur la th\'eorie des immersions. Topology {\bf 1}, 81-100. MR {\bf 27}\,\#2992.
\item{}Po\'enaru, V., (1966) On regular homotopy in codimension 1. Ann.\,of Math. {\bf 83}, 257-265. MR {\bf 33}\,\#732.
\item{}Po\'enaru, V., (1968) Extensions des immersions en codimension 1 (d'apr\`es S. Blank). S\'em. Bourbaki, 
\itemitem{}Expos\'e 342, 1-33. Reprinted in: S\'eminaire Bourbaki, Vol.\,10.\,Ann\'ees 1966/67-1967/68. Soci\'et\'e Math. de France, Paris. (1995), 473-505. MR {\bf 99f:}\,00042.
\item{}Po\'enaru, V., (1971) Homotopy theory and differentiable singularities.\,Manifolds-Amsterdam 1970 
\itemitem{}(Proc.\,Nuffic Summer School) pp.\,106-132.\,Lecture Notes in Math.\,Vol.\,197.\,Springer.\,MR {\bf 44}\,\#2250.
\item{}Smale, S., (1958) A classification of immersions of the two-sphere. Trans. Amer. Math. Soc. {\bf 90}, 
\itemitem{}281-290. MR {\bf 21}\,\#2984.
\item{}Smale, S., (1959) The classification of immersions of spheres in Euclidean spaces. Ann. Math. {\bf 69}, 
\itemitem{}327-344. MR {\bf 21}\,\#3862.
\item{}Spring, D., (1993) Note on the history of immersion theory. From Topology to Computation: Pro-
\itemitem{}ceedings of the Smalefest (M.W. Hirsch, J.E. Marsden, M. Shub, eds.) Springer Verlag, New York, 114-116. MR {\bf 97b:}\,57002.
\item{}Spring, D., (1998) Convex integration theory: solutions to the $h$-principle in geometry and topology. 
\itemitem{}Monographs in Mathematics.\,Vol.\,92. Birkh\"auser Verlag, Basel. MR {\bf 99e:}\,58024. 
\item{}Spring, D., (2002) Directed embeddings and the simplification of singularities. Commun. Contemp. 
\itemitem{}Math. {\bf 4}, 107-144. MR {\bf 2003h:}\,58055.
\item{}Thom, R., (1957) La classification des immersions (d'apr\`es Smale). S\'eminaire Bourbaki, Expos\'e 157,
\itemitem{}11pp. MR {\bf 21}\,\#4429.
\item{}Thom, R., (1959) Remarques sur les probl\`emes comportant des in\'equations diff\'erentielles globales. 
\itemitem{}Bull. Soc. Math. France {\bf 87}, 455-461. MR {\bf 22}\,\#12537.
\item{}Thurston, W.P., (1974) The theory of foliations of codimension greater than one. Comment. Math. 
\itemitem{}Helv. {\bf 49}, 214-231. MR {\bf 51}\,\#6846.
\item{}Thurston, W.P., (1976) Existence of codimension-one foliations. Ann.\,of Math. {\bf 104}, 249-268. 
\itemitem{}MR {\bf 54}\,\#13934.  
\item{}Whitney, H., (1937) On regular closed curves in the plane. Compositio Math. {\bf 4}, 276-284.
\item{}Whitney, H., (1944) The singularities of a smooth $n$-manifold in $(2n-1)$-space. Ann.\,of Math. {\bf 45}, 
\itemitem{}247-293. MR {\bf 5},\,274a.
\item{}Whitney, H., (1955) On singularities of mappings of euclidean spaces.\,I. Mappings of the plane into 
\itemitem{}the plane. Ann.\,of Math. {\bf 62}, 374-410. MR {\bf 17},\,518d. 

\enddocument